\title{A note on the frontier of a branching reflected Brownian motion}
\author{Wenpin Tang \thanks{ Statistics department, University of California, Berkeley. Email: wenpintang@stat.berkeley.edu.}}
\date{\today}
\documentclass[12pt,a4paper]{article}
\usepackage{graphicx} \usepackage{amsmath, amsfonts,amssymb}
\usepackage{float}
\usepackage{enumerate}
\usepackage[utf8]{inputenc} \usepackage[T1]{fontenc} \usepackage[french,english]{babel} \usepackage{color}
\usepackage[top=1in, bottom=1in, left=1in, right=1in]{geometry}
\DeclareMathOperator\erf{erf}
\newtheorem{theorem}{Theorem}[section]
\newtheorem{lemma}[theorem]{Lemma}
\newtheorem{proposition}[theorem]{Proposition}
\newtheorem{corollary}[theorem]{Corollary}
\newtheorem{conjecture}[theorem]{Conjecture}

\begin{document}
\maketitle
\textbf{Abstract:} In this note, we study the asymptotical frontier behavior of a branching reflected Brownian motion. There is essentially no difference in maximal displacement between a branching Brownian motion and its reflected counterpart. We provide two proofs of this fact, one via a soft argument on the dependance of two-sided extremal particles in a branching Brownian motion and the other based on direct computations as in Roberts \cite{Roberts}. The asymptotics of minimal displacement is also given.


\setcounter{tocdepth}{1}
%---------------------------------------------------------------------------------------------------------------------------------------
\section{Introduction and main results}
The study of the frontier behavior of a branching Brownian motion with quadratic branching mechanism was initiated by McKean \cite{Mckean}, where a renewal argument shows that the cumulative distribution of its maximal displacement $M_t$ at time $t$, $u(t,x):=\mathbb{P}(M_t \leq x)$ solves the semilinear heat equation
\begin{equation}
\label{FKPP}
\partial_t u=\frac{1}{2} \partial_{xx} u+u^2-u;
\end{equation} 
with Heaviside initial condition
\begin{equation}
\label{Heaviside}
u(0,x)=\left\{ \begin{array}{rcl}
        1 & \mbox{if}
         & x \geq 0, \\ 0  & \mbox{if} & x<0.
                \end{array}\right.
\end{equation}

The equation \eqref{FKPP} was first carried out by Fisher \cite{Fisher} and then fully investigated by Kolmogorov et al \cite{KPP} (which is known as F-KPP equation). It can be deduced from analytical results that the extremal particle in a branching Brownian motion sits around $\sqrt{2}t$ as $t \rightarrow \infty$. Later Bramson \cite{Bramson2,Bramson1} provided higher order expansions using a probabilistic approach. Hu and Shi \cite{HuShi} proved an almost sure fluctuation result in the branching random walks setting and Roberts \cite{Roberts} derived a similar result related to branching Brownian motion. We summarize their results as follows.
\begin{theorem}\label{summary} Let $(M_t;t \geq 0)$ be the maximal displacement of a branching Brownian motion and $q_{\delta}(t):=\sup\left\{x \geq 0; ~\mathbb{P}(M_t  \leq x)\leq \delta\right\}$ be its $\delta-$quantile for $0<\delta<1$. Then
\begin{enumerate}[(i).]
\item
\cite{Bramson1}
\begin{equation}
\label{morder}
q_{\delta}(t) = \sqrt{2}t-\frac{3}{2\sqrt{2}}\log t +\mathcal{O}(1) \quad \mbox{as}~t \rightarrow \infty,
\end{equation}
\item
\cite{HuShi,Roberts}
\begin{equation}
\label{liminf}
\liminf_{t \rightarrow \infty} \frac{M_t-\sqrt{2}t}{\log t}=-\frac{3}{2\sqrt{2}} \quad a.s.
\end{equation}
and
\begin{equation}
\label{limsup}
\limsup_{t \rightarrow \infty} \frac{M_t-\sqrt{2}t}{\log t}=-\frac{1}{2\sqrt{2}} \quad a.s.
\end{equation}
\end{enumerate}
\end{theorem}

In the note, we deal with a related problem on the asymptotical frontier behavior of a branching reflected Brownian motion. This work can be viewed as a modest complement to Arguin et al \cite{ABK3} and Roberts \cite{Roberts}, from which the strategies borrow mainly with appropriate modifications and refinements.

We consider a continuous-time branching reflected Brownian motion with quadratic mechanism: the system starts with a single particle at the origin and follows a reflected Brownian motion with zero drift and unit variance. After an exponential time with parameter $1$, it splits into two new particles, each of which -- relative to their common birth place -- moves as independent copies of reflected Brownian motion and branches at rate $1$ into two copies of themselves$\dots$etc. The branching system is supercritical and the expected number of particles alive at time $t$, $\mathbb{E} \# N^R(t)=e^t$ for $t \geq 0$.

Denote $X_v^R(t)$ the position of $v \in N^R(t)$ alive at time $t$ in the branching reflected Brownian motion. For $s<t$, $X_v^{R}(s)$ is the position of the ancestor of $v$ that was alive at time $s$. Define $M_t^R:=\underset{v \in N^R(t)}{\max}X^R_v(t)$ its maximal displacement. Again as in Theorem \ref{summary} we are interested in the $\delta$-quantile of the maximal displacement, i.e.
\begin{equation}
\label{median}
q_{\delta}^R(t):=\sup\left\{x \geq 0; ~\mathbb{P}(M^R_t  \leq x)\leq \delta \right\}.
\end{equation}

Observe that the reflected Brownian motion is of non-stationary increments and thus it is not a simple task to identify $u^R(t,x):=\mathbb{P}(M^R_t \leq x)$ as the solution of certain partial differential equation. Nevertheless, it is still possible to associate $u^R(t,x,y):=\mathbb{P}_y(M^R_t \leq x)$ to some partial differential equation, where $\mathbb{P}_y$ is the law of branching reflected Brownian motion starting at $y \geq 0$. We defer the discussion to Section \ref{Open}.

Let us mention some immediate result, at least morally, concerning $(q_{\delta}^{R}(t); t \geq 0)$ and $(M_t^R; t \geq 0)$. Note that the branching reflected Brownian motion can be constructed as the absolute values of all the particles in a branching Brownian motion. Therefore, the rightmost particle in the branching reflected Brownian motion is the maximum of two identically distributed (though not independent) copies of the one-sided extremum. By Theorem \ref{summary},  we see that
\begin{equation*}
q_{\delta}^R(t) \gtrsim \sqrt{2}t-\frac{3}{2\sqrt{2}}\log t +\mathcal{O}(1) \quad \mbox{as}~t \rightarrow \infty.
\end{equation*}
\begin{equation*}
\liminf_{t \rightarrow \infty} \frac{M_t^R-\sqrt{2}t}{\log t} \geq-\frac{3}{2\sqrt{2}} \quad a.s.
\end{equation*}
and
\begin{equation*}
\limsup_{t \rightarrow \infty} \frac{M_t^R-\sqrt{2}t}{\log t}=-\frac{1}{2\sqrt{2}} \quad a.s.
\end{equation*}

Now an interesting question is to determine whether $q_{\delta}^{R}(t)$ and $q_{\delta}(t)$ have exactly the same order as in \eqref{morder} when $t \rightarrow \infty$. The main result of the work provides an affirmative answer to this question.
\begin{theorem}
\label{main}
Let $(q_{\delta}^R(t); t \geq 0)$ defined as in \eqref{median}. Then
\begin{equation}
q_{\delta}^R(t) = \sqrt{2}t-\frac{3}{2\sqrt{2}}\log t +\mathcal{O}(1) \quad \mbox{as}~t \rightarrow \infty.
\end{equation}
\end{theorem}

The above result is not that surprising since in a branching Brownian motion, the rightmost particle is asymptotically independent of the leftmost one. We will develop the circle of ideas in Section \ref{weak} which, together with $(i)$ in Theorem \ref{summary}, gives a proof of Theorem \ref{main}. Remark that the fact is implicitly suggested by Poissonian structure of the extremal process in a branching Brownian motion, see Lalley and Sellke \cite{LS}, Arguin et al \cite{ABK3,ABK1,ABK2} and A\"{i}d\'{e}kon et al \cite{Dabian} for details.

But even without resort to Theorem \ref{summary}, Theorem \ref{main} can still be established using the first/second moment method as in Roberts \cite{Roberts}.  Note that the key idea of Roberts \cite{Roberts} consists in finding the derivative martingale that forces one particle (the so-called spine) to stay below certain curves and by measure change it corresponds to the law of three-dimensional Bessel process. This observation simplifies enormously the computation. However, this method does not work well in the reflected Brownian motion case due to the reflected source -- local times. Thus extra computation is required regarding boundary crossing probability for reflected Brownian motion. All these will be discussed in Section \ref{maximal}. Furthermore, with the ingredients in the proof, we are able to derive an almost sure fluctuation result for $(M^R_t; t \geq 0)$.
\begin{corollary}
\label{corollary}
The maximal displacement of branching reflected Brownian motion $M^R_t$ satisfies
\begin{equation}
\label{liminfr}
\liminf_{t \rightarrow \infty} \frac{M^R_t-\sqrt{2}t}{\log t}=-\frac{3}{2\sqrt{2}} \quad a.s.
\end{equation}
and
\begin{equation}
\label{limsupr}
\limsup_{t \rightarrow \infty} \frac{M^R_t-\sqrt{2}t}{\log t}=-\frac{1}{2\sqrt{2}} \quad a.s.
\end{equation}
\end{corollary}

In a branching Brownian motion, the minimal displacement is opposite to the maximal one and it suffices to study the latter. Nevertheless, in the reflected case, the minimal displacement should be treated separately. The following result regarding the minimal displacement of a branching reflected Brownian motion is a direct consequence of Watanabe's law of large numbers of branching processes \cite{Watanabe}. For the sake of completeness, we include it in Section \ref{minimal}.
\begin{proposition} The minimal displacement of a branching reflected Brownian motion $m^R_t$ satisfies
\begin{equation}
\label{lf}
\lim_{t \rightarrow 0} m^R_t =0 \quad a.s.
\end{equation}
\end{proposition}
%------------------------------------------------------------------------------------------------------------------------------------
\section{Dependance of two-sided extremes in BBM}
\label{weak}
In this section, we provide a proof of Theorem \ref{main} by exploring the dependance of two-sided extremal particles in a branching Brownian motion. We show that at large times, the rightmost extremal particles are asymptotically independent of the leftmost ones. To formulate our result, we need the following notations from Arguin et al \cite{ABK3}.

Let $X_v(t)$ be the position of $v \in N(t)$ alive at time $t$ in the branching Brownian motion. For $s<t$, $X_v(s)$ is the position of the ancestor of $v$ that was alive at time $s$. The correlations among particles at fixed time $t \geq 0$ can be expressed in terms of genealogy distance:
\begin{equation}
\label{Q}
\mathbb{E}[X_u(t)X_v(t)] = Q_t(u,v) \quad \mbox{for}~ u,v \in N(t),
\end{equation}
where $Q_t(u,v):=\sup\{s \leq t; X_u(s)=X_v(s)\} \in [0,t]$ is the most recent common ancestor of $u$ and $v$. Moreover for $a \in \mathbb{R}$ and $f: \mathbb{R}^{+} \rightarrow \mathbb{R}^{+}$ such that $f(t)=o(t)$ as $t \rightarrow \infty$, define
\begin{equation}
\label{cluster}
N^a_f(t):=\{u \in N(t); X_u(t) \in [at-f(t),at+f(t)] \}
\end{equation}
the set of particles falling into the cluster ranging from $at-f(t)$ to $at+f(t)$.

According to Theorem \ref{summary}, the set $N^a_f(t)$ is non-empty if $a \in (-\sqrt{2},\sqrt{2})$ or $a=\pm \sqrt{2}$ and $f(t) \gtrsim \frac{1}{2 \sqrt{2}} \log t$. In the sequel, we suppose that this condition is always satisfied. To abbreviate the notations, we write $X_u(t) \in at+o(t)$ instead of $X_u(t) \in [at-f(t),at+f(t)]$ for some valid function $f$. Similarly, we denote $N^a$ for $N^a_f$. The following result suggests the genealogy of particles lying in the clusters with different index $a$.
\begin{proposition}
\label{indep}
For $a <b$ and $r>0$ such that $r+\frac{(b-a)^2}{4(1-r)}>2$, we have
\begin{equation*}
\mathbb{P}\Bigg(\exists (u,v) \in N^a(t) \times N^b(t)~\mbox{such that}~Q_t(u,v) \in [rt,t]\Bigg) \rightarrow 0 \quad \mbox{as}~t \rightarrow \infty,
\end{equation*}
where $Q_t$ is defined as in \eqref{Q} and $N^a(t)$, $N^b(t)$ are defined as in \eqref{cluster}.
\end{proposition}
\textbf{Proof:} Using the first moment method, we have
\begin{multline*}
\mathbb{P}\Bigg(\exists (u,v) \in N^a(t) \times N^b(t)~\mbox{such that}~Q_t(u,v) \in [rt,t]\Bigg) \\ \leq \mathbb{E}\Bigg[\#\Bigg\{(u,v) \in N^a(t) \times N^b(t)~\mbox{such that}~Q_t(u,v) \in [rt,t]\Bigg\}\Bigg].
\end{multline*}
According to many-to-two principle (see Lemma $10$ in Bramson \cite{Bramson1} or formulas $(4.3)-(4.4)$ in Arguin et al \cite{ABK3}), 
 \begin{multline}
 \label{many}
\mathbb{E}\Bigg[\#\Bigg\{(u,v) \in N^a(t) \times N^b(t)~\mbox{such that}~Q_t(u,v) \in [rt,t]\Bigg\}\Bigg] \\ =2 e^{2t} \int_{rt}^t e^{-s} ds\int_{\mathbb{R}} \mu_s(dx) \mathbb{P}(X_u(t) \in N^a(t)|X_u(s)=x)\\ \times \mathbb{P}(X_v(t) \in N^b(t)|X_v(s)=x).
\end{multline}
where $\mu_s(dx)$ is the Gaussian measure with variance $s$.
Note that from the branching point $(s,x)$ of $u$ and $v$, two particles perform independent Brownian motions starting at $x$. Fix $\epsilon>0$ small enough, 
\begin{align}
&\quad ~  \mathbb{P}(X_u(t) \in N^a(t)|X_u(s)=x)  \times \mathbb{P}(X_v(t) \in N^b(t)|X_v(s)=x) \notag\\
&\leq  \mathbb{P}\Bigg(B_0~\mbox{hits}~at+o(t)-x~\mbox{in time}~t-s\Bigg) \times \mathbb{P}\Bigg(B_0~\mbox{hits}~bt+o(t)-x~\mbox{in time}~t-s\Bigg)                \notag\\
&\leq \left\{ \begin{array}{lcl}
         \mathbb{P}\Bigg(|B_0(t-s)| \geq (b-a-\epsilon)t+o(t)\Bigg) \quad  \quad  \quad \quad \quad \quad  \mbox{if}~
         x<(a+\epsilon)t~\mbox{or}~x>(b-\epsilon)t \\ \mathbb{P}\Bigg(|B_0(t-s)| \geq x-at+o(t)\Bigg) \times \mathbb{P}\Bigg(|B_0(t-s)| \geq bt-x+o(t)\Bigg)  ~~~~~\mbox{otherwise}.               
\end{array}\right. \notag
\end{align}
where $B_0$ in the above expressions is a standard Brownian motion starting at $0$. It is classical that for $z \geq 0$, $\mathbb{P}(\mathcal{N}(0,1) \geq z) \leq \frac{1}{\sqrt{2 \pi z}} \exp(-\frac{z^2}{2})$ where $\mathcal{N}(0,1)$ is standard normal distribution. Therefore for $rt \leq s \leq t$,
\begin{align}
\mathbb{P}\Bigg(\left|B_0(t-s)\right| \geq (b-a)t+o(t) \Bigg)&= \mathbb{P} \left( \mathcal{N}(0,1) \geq \frac{(b-a-\epsilon)t+o(t)}{\sqrt{t-s}}\right) \notag\\
& \label{12}\leq \frac{1}{\sqrt{2\pi(b-a-\epsilon)}}t^{-\frac{1}{4}}  \exp\left(-\frac{(b-a-\epsilon)^2t}{2(1-r)}+o(t)\right)
\end{align}
and for $(a+\epsilon)t \leq x \leq (b-\epsilon)t$,
\begin{align}
& \quad~ \mathbb{P}\Bigg(|B_0(t-s)| \geq x-at+o(t)\Bigg) \times \mathbb{P}\Bigg(|B_0(t-s)| \geq bt-x+o(t)\Bigg) \notag\\
&= \mathbb{P} \left( \mathcal{N}(0,1) \geq \frac{x-at+o(t)}{\sqrt{t-s}}\right) \times \mathbb{P} \left( \mathcal{N}(0,1) \geq \frac{bt-x+o(t)}{\sqrt{t-s}}\right)\notag\\
& \leq \frac{1}{2 \pi \epsilon}t^{-\frac{1}{2}}  \exp\left(-\frac{(x-at+o(t))^2+(bt-x+o(t))^2}{2(1-r)t}\right) \notag\\
& \label{13} \leq \frac{1}{2 \pi \epsilon}t^{-\frac{1}{2}}  \exp\left(-\frac{(b-a)^2t}{4(1-r)}+o(t)\right)
\end{align}
Injecting \eqref{12} and \eqref{13} into \eqref{many}, we obtain
\begin{align*}
& \quad ~\mathbb{E}\Bigg[\#\Bigg\{(u,v) \in N^a(t) \times N^b(t)~\mbox{such that}~Q_t(u,v) \in [rt,t]\Bigg\}\Bigg] \\
&\leq K t^{-\frac{1}{4}} \exp\left[\left(2-r-\frac{(b-a)^2}{4(1-r)}\right)t+o(t)\right] \rightarrow 0 \quad \mbox{as}~t \rightarrow \infty. \quad \square
\end{align*}

According to $(ii)$ in Theorem \ref{summary}, the maximal displacement $M(t)$ and the minimal displacement $m(t)$) in a branching Brownian motion satisfy
$$M(t)=\sqrt{2}t+o(t) \quad \mbox{and} \quad m(t)=-\sqrt{2}t+o(t).$$
Thus, $M(t) \in N^{\sqrt{2}}(t)$ and $m(t) \in N^{-\sqrt{2}}(t)$. We check that $r+\frac{2}{1-r}>2$ for all $r>0$. Proposition \ref{indep} then leads to asymptotic independence of these two particles.
\begin{corollary} \label{keycor} For $r>0$, $M(t)$ (resp. $m(t)$) the maximal displacement (resp. the minimal displacement) in a branching Brownian motion, we have
\begin{equation*}
Q_t(M(t),m(t)) \leq rt \quad \mbox{for}~t~\mbox{large enough},
\end{equation*}
where $Q_t$ is defined as in \eqref{Q}. In other words, for $x>0$,
\begin{equation*}
\mathbb{P}\Bigg(M_t \leq x~\mbox{and}~m_t \geq -x\Bigg)-\mathbb{P}(M_t \leq x)^2 \rightarrow 0 \quad \mbox{as}~t \rightarrow \infty,
\end{equation*}
\end{corollary}
\textbf{First proof of Theorem \ref{main}:} It is already clear from the construction of a branching reflected Brownian motion that $q^R_{\delta}(t) \geq q_{\delta}(t)$. Suppose now $\mathbb{P}(M^R_t \leq x) \leq \delta$. According to Corollary \ref{keycor}, for arbitrary small $\epsilon>0$, $$\mathbb{P}(M_t \leq x)^2 -\mathbb{P}(M_t^R \leq x) \leq \epsilon \quad \mbox{for}~t~\mbox{large enough}.$$
Consequently, $\mathbb{P}(M_t \leq x) \leq \sqrt{\delta+\epsilon}$ and thus $q^R_{\delta}(t) \leq q_{\sqrt{\delta + \epsilon}}(t)$ when $t$ is large. Then $(i)$ in Theorem \ref{summary} permits to conclude.  $\square$
%------------------------------------------------------------------------------------------------------------------------------------
\section{Maximal displacement in BRBM}
\label{maximal}
In this part, we provide an alternative approach to Theorem \ref{main} as well as Corollary \ref{corollary} without appealing to Theorem \ref{summary}. The main argument is the same as in Roberts \cite{Roberts} but extra efforts are needed as explained in the introduction. We hope that the method is still applicable when dealing with the frontier behavior of some other branching diffusions, e.g. branching $\delta-$Bessel processes (see Chapter XI, Revuz and Yor \cite{RY} for backgrounds).
\subsection{Background and basic tools}
\label{premin}
%-------------------------------------------------------------------------------------
\subsubsection{Reflected Brownian motion}
We recall some basic properties of reflected Brownian motion. A reflected Brownian motion $(X_t;t \geq 0)$ is defined as the unique strong solution of the Skorokhod equation
$$X_t=B_t+L_t,$$
where $(B_t;t \geq 0)$ is standard Brownian motion and $(L_t;t \geq 0)$ is local times process, which is increasing and whose measure is supported on the zero set of $(X_t;t \geq 0)$. In particular, a reflected Brownian motion has the same law as the absolute value of linear Brownian motion. By abuse of language, we identify the reflected Brownian motion $(X_t;t \geq 0)$ with $(|B_t|;t \geq 0)$ in the sequel.

It is well-known that $(|B_t|;t \geq 0)$ is a strong Markov process with transition density 
\begin{equation}
\label{transden}
p^R(s,x;t,y):=\frac{1}{\sqrt{2 \pi (t-s)}} \left[ \exp\left(-\frac{(y-x)^2}{2(t-s)} \right) + \exp \left(-\frac{(y+x)^2}{2(t-s)} \right) \right],
\end{equation}
for $0 \leq s \leq t$ and $x,y \geq 0$.

As mentioned in the introduction, it is indispensable to compute the probability for a reflected Brownian motion to stay below certain curves and this does not seem to be easily derived by change of measures. More generally, few closed-form expressions are available regarding the general boundary crossing probability for linear Brownian motion. 

Nevertheless in our case, we only need to know some affine boundary crossing probability for reflected Brownian bridges, which has been calculated by Abundo \cite{Abundo}. We also refer readers to Salminen and Yor \cite{SYor} for related results. 
\begin{theorem}  \cite{Abundo}
\label{Abundo}
Let $a,b>0$ and $|x| <a+bt$ for some $t \geq 0$. Define $\tau_{a,b}:=\inf\{s \geq 0; |B_s| \geq a+bs\}$ the first hitting time of reflected Brownian motion to affine boundary. Then
\begin{multline}
\label{boundary}
\mathbb{P}(\tau_{a,b} \geq t~\mbox{and}~|B_t| \in dx)/dx\\=\sqrt{\frac{2}{\pi t}}\exp\left(-\frac{x^2}{2t}\right)\sum_{n \in \mathbb{Z}}(-1)^n \exp\left[-2a\left(b+\frac{a}{t}\right)n^2\right] \cosh\left(\frac{2ax}{t}n\right).
\end{multline}
\end{theorem}
%---------------------------------------------------------------------------------------------
\subsubsection{Frontier of branching Brownian motion}
We summarize some previous results, especially some key estimations regarding the frontier of branching Brownian motion. All these can be read from Roberts \cite{Roberts}, which provides a much simpler approach compared to Bramson's original proof \cite{Bramson1} of $(i)$ in Theorem \ref{summary}.

Fix $t>0$ and $y \in \mathbb{R}$, the following two sets are of special importance in the proof:
\begin{equation}
\label{h}
H(y,t):=\#\{u \in N(t); X_u(s) \leq \beta s+1~\forall s \leq t~\mbox{and}~\beta t-1 \leq X_u(t) \leq \beta t\},
\end{equation}
and
\begin{equation}
\label{gamma}
\Gamma(y,t):=\#\{u \in N(t); X_u(s) \leq \beta s+L(s)+y+1~\forall s \leq t~\mbox{and}~\beta t-1 \leq X_u(t) \leq \beta t+y\},
\end{equation}
where 
\begin{equation}
\label{beta}
\beta:=\sqrt{2}-\frac{3}{2 \sqrt{2}} \frac{\log t}{t}+\frac{y}{t},
\end{equation}
and $L \in \mathcal{C}^{\infty}(\mathbb{R})$ satisfying
\begin{equation}
\label{L}
L(s):=\left\{ \begin{array}{ccl}
        \frac{3}{2 \sqrt{2}}\log(s+1) & \mbox{for}
         & s \in [0,\frac{t}{2}-1], \\ \frac{3}{2 \sqrt{2}}\log(t-s+1) & \mbox{for} & s \in [\frac{t}{2}+1,t],
                \end{array}\right.
\end{equation}
with $L^{''}(s) \in [-\frac{10}{t},0]$ for $s \in [\frac{t}{2}-1,\frac{t}{2}+1]$.

On one hand, using the second moment method for $H(y,t)$, it can be proved that there exists $C_H>0$ such that $\mathbb{P}(H(y,t)>0) \geq C_He^{-\sqrt{2}y}$. Hence,
\begin{equation}
\label{LowerBBM}
\mathbb{P}\left(M_t \geq \sqrt{2}t-\frac{3}{2\sqrt{2}}\log t +y\right) \geq C_He^{-\sqrt{2}y}.
\end{equation}
And on the other hand, the first moment method for $\Gamma(y,t)$ together with some coupling arguments guarantee the existence of $C_{\Gamma}>0$ such that 
\begin{equation}
\label{UpperBBM}
\mathbb{P}\left(M_t \geq \sqrt{2}t-\frac{3}{2\sqrt{2}}\log t +y\right) \leq C_{\Gamma} (y+2)^2 e^{-\sqrt{2}y}.
\end{equation}
Then $(i)$ of Theorem \ref{summary} follows directly \eqref{LowerBBM} and \eqref{UpperBBM}. With some extra efforts, the almost sure fluctuation result of $M_t$, i.e. $(ii)$ of Theorem \ref{summary} is also derived via these estimations.
%-------------------------------------------------------------------------------------------------
\subsection{Asymptotics of $q_{\delta}(t)$ -- Proof of Theorem \ref{main}}
The current part is devoted to the proof of our main result, Theorem \ref{main}. As in the branching Brownian motion case, we subdivide the proof into two propositions.
\begin{proposition}
\label{Lower}
There exists $C^R_H>0$ such that for $t \geq 1$ and $y \in [0,\sqrt{t}]$,
\begin{equation} 
\label{LowerBRBM}
\mathbb{P}\left(M^R_t \geq \sqrt{2}t-\frac{3}{2\sqrt{2}}\log t +y\right) \geq C^R_He^{-\sqrt{2}y}.
\end{equation}
\end{proposition}
\begin{proposition}
\label{Upper}
There exists $C^R_{\Gamma}>0$ such that for $t \geq 1$ and $y \in [0,\sqrt{t}]$,
\begin{equation}
\label{UpperBRBM}
\mathbb{P}\left(M^R_t \geq \sqrt{2}t-\frac{3}{2\sqrt{2}}\log t +y\right) \leq C^R_{\Gamma} (y+2)^2 e^{-\sqrt{2}y}.
\end{equation}
\end{proposition}
Before proving these two results, let us indicate how we use them to prove Theorem \ref{main}.\\\\
\textbf{Second proof of Theorem \ref{main}:} From \eqref{LowerBRBM} and \eqref{UpperBRBM} follows for $t \geq 1$ and $y \in[0,\sqrt{t}]$,
\begin{equation}
\label{twoside}
1-C_{\Gamma}^R(y+2)^2\ e^{-\sqrt{2}y} \leq \mathbb{P}\left(M^R_t \leq \sqrt{2}t-\frac{3}{2 \sqrt{2}}\log t+y\right) \leq 1-C_H^R e^{-\sqrt{2}y}.
\end{equation}
Take $y \sim \sqrt{t}$ and $t \rightarrow \infty$, the bounds in both sides of \eqref{twoside} converge to $1$. Thus there exists $\delta_0>0$ such that for $\delta_0 \leq \delta<1$,
$$q_{\delta}^R(t)=\sqrt{2}t-\frac{3}{2 \sqrt{2}}\log t+\mathcal{O}(1).$$
In addition, $q_{\delta}^R(t) \leq q_{\delta_0}^R(t)=\sqrt{2}t-\frac{3}{2 \sqrt{2}}\log t+\mathcal{O}(1)$ for $0 \leq \delta<\delta_0$. 

Now fix $\delta>0$ and $\epsilon>0$. Choose $L>0$ such that $\mathbb{E}(\delta^{\#N^{R}(L)}) \leq \frac{\epsilon}{2}$ and $a>0$ such that $\mathbb{P}(M_L^R>a) \leq \frac{\epsilon}{2}$. We have for $t \geq L$,
\begin{align}
&\quad~ \mathbb{P}(M_t^R \geq q_{\delta}^R(t-L)+a) \notag\\ & =\mathbb{P}\left(\max_{u \in N^R(L)} \max_{v \leftarrow u}X_v^R(t) \geq q_{\delta}^R(t-L)+a\right) \notag\\
                                                                           & =\mathbb{P}(M_L^R>a) +\mathbb{P}\left(M^R_L \leq a~\mbox{and}~\max_{u \in N^R(L)} \max_{v \leftarrow u}X_v^R(t) \geq q_{\delta}^R(t-L)+a\right) \notag\\
                                                                           & \leq \frac{\epsilon}{2} + \mathbb{E}[\mathbb{P}(M^R_{t-L}>q_{\delta}^R(t-L))^{\#N^R(L)}] \leq \epsilon, \label{tight1}
\end{align}
where $v \leftarrow u$ means that $v$ is a descendent of $u$. Similarly, we can choose $\tilde{L}$ such that
\begin{equation}
\label{tight2}
P(M_t^R \leq q_{\delta}(t-\tilde{L})) \leq \epsilon.
\end{equation}
By \eqref{tight1} and \eqref{tight2}, $(M_t^R-q_{\delta}^R(t); t \geq 0)$ is tight for $\delta \in (0,1)$. As a consequence, $q_{\delta_0}^R(t)-q_{\delta}^R(t)=\mathcal{O}(1)$ for $0 \leq \delta<\delta_0$, which proves the desired result.  $\square$

The rest of the section is devoted to the proofs of Proposition \ref{Lower} (Section \ref{31}), Proposition \ref{Upper} (Section \ref{32}) and Corollary \ref{corollary} (Section \ref{33}).
%----------------------------------------------------------
\subsubsection{Lower bound -- Proof of Proposition \ref{Lower}}
\label{31}
The proof of Proposition \ref{Lower} is quite simple. In fact, we need to show that
\begin{lemma}
\label{comparison}
For $t \geq 0$, $M^R_t$ is stochastically larger than $M_t$, i.e. for all $x \geq 0$,
$$\mathbb{P}(M_t^R \geq x) \geq \mathbb{P}(M_t \geq x).$$
\end{lemma}
Then \eqref{LowerBRBM} follows immediately \eqref{LowerBBM} with $C^R_H=C_H$ by taking $x=\sqrt{2}t-\frac{3}{2 \sqrt{2}} \log t+y$.\\\\
\textbf{Proof:} Note that $M^R_t \stackrel{d}{=} \max\{M_t,-m_t\}$ where $M_t$ (resp. $m_t$) is the maximal displacement (resp. minimal displacement) in a branching Brownian motion.  $\square$
%--------------------------------------------------------
\subsubsection{Upper bound -- Proof of Proposition \ref{Upper}}
\label{32}
We derive the upper bound for the maximal displacement of a branching reflected Brownian motion, i.e. Proposition \ref{Upper}. 

Recall that the sets $H^{R}(y,t)$ and $\Gamma^R(y,t)$ are defined as in \eqref{h} and \eqref{gamma}, in which particles are governed by reflected Brownian motion $X^R$ instead of standard Brownian motion $X$. We need some estimations for $H^R$ and $\Gamma^R$.
\begin{lemma}
\label{himp}
There exists $c_1,C_1 > 0$ such that for $t \geq 1$ and $y \in [0,\sqrt{t}]$,
\begin{equation}
\label{esth}
c_1e^{-\sqrt{2}y} \leq \mathbb{E}H^R(y,t) \leq C_1 e^{-\sqrt{2}y}.
\end{equation}
\end{lemma}
\begin{lemma}
\label{gammaimp}
There exists $C_2 > 0$ such that for $t \geq 1$ and $y \in [0,\sqrt{t}]$,
\label{estimgamma}
\begin{equation}
\label{estgamma}
\mathbb{E}\Gamma^R(y,t) \leq C_2(y+2)^2 e^{-\sqrt{2}y}.
\end{equation}
\end{lemma}
\textbf{Proof of Proposition \ref{Upper}:} \eqref{UpperBRBM} follows \eqref{esth} and \eqref{estgamma} in the same way as the proof of Proposition $11$, Roberts \cite{Roberts}.  $\square$\\\\
 We now turn to prove Lemma \ref{himp} and Lemma \ref{gammaimp}.\\\\
 \textbf{Proof of Lemma \ref{himp}:} According to many-to-one principle,
 \begin{equation}
 \label{expectationh}
 \mathbb{E}H^R(y,t)=e^t \mathbb{P}(\tau_{1,\beta} \geq t~\mbox{and}~|B_t| \in [\beta t-1,\beta t]),
 \end{equation}
 where $\tau_{1,\beta}$ is defined as in Theorem \ref{Abundo}. By \eqref{boundary},
 \begin{multline}
 \label{boundarybis}
\mathbb{P}(\tau_{1,\beta} \geq t~\mbox{and}~|B_t| \in [\beta t-1,\beta t])\\ =\sqrt{\frac{2}{\pi t}}\sum_{n \in \mathbb{Z}} (-1)^n \exp\left[-2\left(\beta+\frac{1}{t}\right)n^2\right] \int_{\beta t-1}^{\beta t} \cosh \left(\frac{2nx}{t}\right) \exp\left(-\frac{x^2}{2t}\right)
 \end{multline}
 Now the key issue is to evaluate the asymptotical behavior of the integral in \eqref{boundarybis}.
 \begin{multline}
 \label{310}
 \int_{\beta t-1}^{\beta t} \cosh \left(\frac{2nx}{t}\right) \exp\left(-\frac{x^2}{2t}\right) \\ =  \sqrt{\frac{\pi t}{8}} \exp\left(\frac{2n^2}{t}\right) \Bigg[ -\erf \left(\beta \sqrt{\frac{t}{2}}+n \sqrt{\frac{2}{t}} -\sqrt{\frac{1}{2t}} \right)-\erf \left(\beta \sqrt{\frac{t}{2}}-n \sqrt{\frac{2}{t}}-\sqrt{\frac{1}{2t}} \right) \\+ \erf \left(\beta \sqrt{\frac{t}{2}}+n \sqrt{\frac{2}{t}}\right) + \erf \left(\beta \sqrt{\frac{t}{2}}-n \sqrt{\frac{2}{t}}\right)\Bigg],
 \end{multline}
 whereas
 \begin{equation}
 \label{320}
\erf \left(\beta \sqrt{\frac{t}{2}}+n \sqrt{\frac{2}{t}}\right) + \erf \left(\beta \sqrt{\frac{t}{2}}-n \sqrt{\frac{2}{t}}\right) = 2-\frac{2}{\sqrt{\pi}} e^{-t}(t+n^2) \cosh(\sqrt{8} n)+o(e^{-t}),
\end{equation}
and 
\begin{multline}
\label{330}
\erf \left(\beta \sqrt{\frac{t}{2}}+n \sqrt{\frac{2}{t}} -\sqrt{\frac{1}{2t}} \right)+\erf \left(\beta \sqrt{\frac{t}{2}}-n \sqrt{\frac{2}{t}}-\sqrt{\frac{1}{2t}} \right)\\=2-\frac{2}{\sqrt{\pi}}e^{\sqrt{2}-t}(t+n^2)\cosh(\sqrt{8} n)+o(e^{-t}).
\end{multline}
Injecting \eqref{310}, \eqref{320} and \eqref{330} into \eqref{boundarybis}, we obtain 
\begin{multline}
\label{340}
\mathbb{P}(\tau_{1,\beta} \geq t~\mbox{and}~|B_t| \in [\beta t-1,\beta t]) \\=\frac{e^{\sqrt{2}}-1}{\sqrt{\pi}}e^{-t}\Bigg[t \sum_{n \in \mathbb{Z}}(-1)^n e^{-\sqrt{8}n^2} \cosh(\sqrt{8}n) \\+ \sum_{n \in \mathbb{Z}} (-1)^n n^2 e^{-\sqrt{8}n^2} \cosh(\sqrt{8}n)+o(1)\Bigg] \asymp e^{-t},
\end{multline}
since $\underset{n \in \mathbb{Z}}{\sum}(-1)^n e^{-\sqrt{8}n^2} \cosh(\sqrt{8}n)=0$ and $\underset{n \in \mathbb{Z}}{\sum}(-1)^n n^2 e^{-\sqrt{8}n^2} \cosh(\sqrt{8}n) \in ]0, \infty[.$ Then \eqref{esth} follows immediately \eqref{expectationh} and \eqref{340}.  $\square$ \\\\
 \textbf{Proof of Lemme \ref{gammaimp}:} Again by many-to-one principle, we have
 \begin{align}
 \mathbb{E}\Gamma^R(y,t) &=e^t \mathbb{P}(|B_s| \leq \beta s +L(s)+y+1~\forall s \leq t~\mbox{and}~ |B_t| \in[\beta t-1 , \beta t+y]) \notag\\
                                                 & \leq e^t \Bigg[\mathbb{P}(B_s \leq \beta s +L(s)+y+1~\forall s \leq t~\mbox{and}~B_t \in[\beta t-1 , \beta t+y]) \notag\\ 
                                                 &\quad \quad~ +\mathbb{P}(B_s \geq -\beta s -L(s)-y-1~\forall s \leq t~\mbox{and}~ B_s \in [-\beta t-y , -\beta t+1]) \Bigg] \notag\\
                                                 &=2\mathbb{E}\Gamma(y,t) \leq 2C_{\Gamma}(y+2)^2 e^{-\sqrt{2}y}, \notag
 \end{align}
 where the last inequality is due to \eqref{UpperBBM}. It suffices to take $C_2=2 C_{\Gamma}$ in \eqref{estgamma}.  $\square$
 %---------------------------------------------------------
\subsection{Almost sure fluctuation -- Proof of Corollary \ref{corollary}}
\label{33}
\textbf{Proof of Corollary \ref{corollary}:} Following the stochastic comparison in Lemma 
\ref{comparison} together with \eqref{liminf} and \eqref{limsup} in Theorem \ref{summary}, it is straightforward that
\begin{lemma}
\label{360}
\begin{equation*}
\liminf_{t \rightarrow \infty} \frac{M^R_t-\sqrt{2}t}{\log t} \geq -\frac{3}{2\sqrt{2}} \quad a.s.
\end{equation*}
and
\begin{equation*}
\limsup_{t \rightarrow \infty} \frac{M^R_t-\sqrt{2}t}{\log t} \geq -\frac{1}{2\sqrt{2}} \quad a.s.
\end{equation*}
\end{lemma}
Note in addition that the upper bounds for the almost sure fluctuation can be derived in the same way as in Lemma $12$ and Lemma $13$, Roberts \cite{Roberts}. In fact, \eqref{UpperBBM} is the only ingredient that was used to obtain these bounds, which is replaced in our case by \eqref{UpperBRBM} in Proposition \ref{Upper}. Thus,
\begin{lemma}
\label{370}
\begin{equation*}
\liminf_{t \rightarrow \infty} \frac{M^R_t-\sqrt{2}t}{\log t} \leq -\frac{3}{2\sqrt{2}} \quad a.s.
\end{equation*}
and
\begin{equation*}
\limsup_{t \rightarrow \infty} \frac{M^R_t-\sqrt{2}t}{\log t} \leq -\frac{1}{2\sqrt{2}} \quad a.s.
\end{equation*}
\end{lemma}
From Lemma \ref{360} and Lemma \ref{370} follows Corollary \ref{corollary}.  $\square$ 
%------------------------------------------------------------------------------------------------------------------------------------
\section{Minimal displacement in BRBM}
\label{minimal}
We include in this section Watanabe's law of large numbers of branching processes which leads naturally to Proposition \ref{minimal}.

Recall that $N(t)$ is the number of particles alive at time $t$ in a quadratic branching Brownian motion. It is classical that $(e^{-t}N(t); t \geq 0)$ is a positive martingale. By martingale convergence theorem, there exists a random variable $W_{\infty} \in [0, \infty)$ such that
\begin{equation}
\label{W}
e^{-t}N(t) \rightarrow W_{\infty} \quad a.s. \quad \mbox{as}~t \rightarrow \infty.
\end{equation}
In addition, $\mathbb{P}(W_{\infty}=0)$ is the smallest root of $\Phi(s)=s$ in $[0,1]$, where $\Phi(s)=s^2$ for quadratic branching mechanism. Thus, $W_{\infty}>0$ a.s. 

The following theorem due to Watanabe studied the asymptotical behavior of the number of particles confined in any region at large times.
\begin{theorem} \cite{Watanabe}
\label{Watanabe}
For $D \subset \mathbb{R}$, $N_D(t)$ denote the number of particles lying in domain $D$ at time $t$  in a quadratic branching Brownian motion. Then
$$\frac{N_D(t)}{e^t t^{-\frac{1}{2}}} \rightarrow \frac{|D|}{\sqrt{2 \pi}} W_{\infty} \quad a.s. \quad \mbox{as}~t \rightarrow \infty,$$
where $|D|$ is the Lebesgue measure of domain $D$ and $W_{\infty}>0$ a.s. is defined as in \eqref{W}.
\end{theorem}
\textbf{Proof of Proposition \ref{minimal}:} It suffices to take $D^{\epsilon}:=[-\epsilon, \epsilon]$. By Theorem \ref{Watanabe},
$$\frac{N_{D^{\epsilon}}(t)}{e^t t^{-\frac{1}{2}}} \rightarrow \frac{|D|}{\sqrt{2 \pi}} W_{\infty}>0 \quad a.s. \quad \mbox{as}~t \rightarrow \infty.$$
Thus for arbitrary small $\epsilon>0$, $N_{D^{\epsilon}}(t) \neq 0$ for $t$ large enough and $\underset{t \rightarrow \infty}{\limsup} m^R_t \leq \epsilon$.  $\square$
 %------------------------------------------------------------------------------------------------------------------------------------
\section{Appendix: BRBM and PDEs}
\label{Open}
As explained in the introduction, a powerful tool to study the frontier behavior of a branching Brownian motion is the F-KPP equation \eqref{FKPP} with Heaviside initial condition \eqref{Heaviside}, to which $u(t,x):=\mathbb{P}(M_t \leq x)$ is solution. This provides an analytical way to tackle down the problem. 

It is natural to ask whether it is also possible to find certain PDE such that $u^R(t,x):=\mathbb{P}(M^R_t \leq x)$ is solution. The question is subtle since the reflected Brownian motion is space-inhomogeneous in the sense that for some $t_0>0$, $(|B_{t_0+t}|-|B_{t_0}|; t \geq 0)$ is no longer a reflected Brownian motion. As a result, it might be a bad strategy without precising the starting position of branching reflected Brownian motion. Of course, we do not exempt the possibility to find some eligible PDE; however it is far more than obvious.

Now let us take the starting point of reflected Brownian motion into consideration. Define $u^R(t,x,y):=\mathbb{P}_y(M^R_t \leq x)$, where $\mathbb{P}_y$ is the law of branching reflected Brownian motion starting at $y \geq 0$. It is not difficult to adapt the renewal argument in McKean \cite{Mckean} to see that $u^R(t,x,y)$ solves the integral equation
\begin{multline}
\label{integral}
u^R(t,x,y)=e^{-t} \int_{-x}^{\infty} p^R(0,x_0;t,z+x) H(z)dz\\
+e^{-t} \int_0^t \int_0^{\infty} e^{t'}p^R(t',x_0;t,y) \Bigg[u^R(t',z,x)\Bigg]^2 dzdt',
\end{multline}
where $p^R$ is the density of reflected Brownian motion defined as in \eqref{transden} and $H$ is Heaviside function defined as in \eqref{Heaviside}. A direct computation confirms that $\mathbb{R}^{+} \times \mathbb{R}^{+} \times \mathbb{R}^{+} \ni (t,x,y) \rightarrow u^R(t,x,y) \in [0,1]$ is solution to
\begin{equation}
\label{FKPPbis}
\partial_t u^R =\frac{1}{2} \partial_{yy} u^R+(u^R)^2-u^R;
\end{equation}
with Heavide-type initial condition
\begin{equation}
\label{Heavisebis}
u^{R}(0,x,y)=H(x-y);
\end{equation}
and Neumann boundary condition
 \begin{equation}
\label{Neumann}
\partial_{y}u^{R}|_{y=0}=0.
\end{equation}
The initial-boundary PDE system \eqref{FKPPbis}, \eqref{Heavisebis} and \eqref{Neumann} looks similar to the standard F-KPP equation \eqref{FKPP} and \eqref{Heaviside} but is totally different. The following conjecture appears to be reasonable in view of Theorem \ref{main} in the current work.
\begin{conjecture}
\label{conj1}
Let $y \geq 0$ and $q_{\delta}^R(y,t)$ defined as in \eqref{median} starting at $y$. Then
\begin{equation}
u^R(t,x+q^R_{\delta}(y,t),y) \rightarrow w^R(x),
\end{equation}
where $w^R: \mathbb{R}^{+} \rightarrow [0,1]$ does not depend on $y$ and is the solution to certain ODE to be precised.
\end{conjecture}
$$$$
\textbf{Acknowledgement:} The author would like to thank Jim Pitman and Matthew Roberts for enlightening discussions on the subject.
%------------------------------------------------------------------------------------------------------------------------------------
\addcontentsline{toc}{section}{References}
\bibliographystyle{plain}
\bibliography{BRBM}
\end{document}